\documentclass[a4paper]{article}
\usepackage{amsfonts, latexsym, amsthm}

\newcommand{\ord}{\mathop{\rm ord}\nolimits}
\newcommand{\Image}{\mathop{\rm Im}\nolimits}
\newcommand{\scr}{\scriptstyle}

\newtheorem*{theorem}{Theorem}
\newtheorem{lemma}{Lemma}

\newenvironment{definition}
{\smallskip\noindent{\bf Definition\/}:}{\smallskip\par}
\title{Simple singularities of multigerms of curves \thanks{Accepted
       for publication in Revista Matematica Complutense}}
\author{P.~A.~Kolgushkin \and R.~R.~Sadykov}
\date{}

\begin{document}
\maketitle

\begin{abstract}
We classify stably simple reducible curve singularities in complex
spaces of any dimension. This extends the same classification of
of irreducible curve singularities obtained by V.~I.~Arnold. The
proof is essentially based on the method of complete transversals
by J.~Bruce et al.
\end{abstract}

\section {Introduction}
Classification of simple curve singularities has been discussed in
a number of papers. J.~W.~Bruce and T.~J.~Gaffney have classified
irreducible plane curves in~\cite{BG}. In~\cite{GH} C.~G.~Gibson and
C.~A.~Hobbs gave the classification of irreducible simple curve
singularities in the 3-dimensional complex space. M.~Giusti
in~\cite{Giu} classified simple complete intersection singularities in
the $3$-dimensional complex space. In~\cite{Frh} A.~Fr\"uhbis-Krger
classified so called simple determinantal singularities.
Classification of
irreducible (stably) simple curve singularities in a linear complex
space of any dimension was made by V.~I.~Arnold~\cite{Arn}.
We consider reducible curve singularities in a linear complex space of
any dimension and give the list of stably simple ones.

An irreducible curve singularity at the origin in ${\mathbb C}^n$ can
be described by a germ of a complex analytic map
$f: ({\mathbb C}, 0)\to ({\mathbb C}^n, 0)$. Let $L$ (respectively $R$)
be the group of co-ordinate changes in $({\mathbb C}^n, 0)$
(respectively in $({\mathbb C}, 0)$), i.e., the group of germs of
non-degenerate analytic maps $({\mathbb C}^n, 0)\to ({\mathbb C}^n, 0)$
(respectively $({\mathbb C}, 0)\to ({\mathbb C}, 0)$), let $A=L\times R$.
The group $L$ (respectively $R$) is called  the group of left
(respectively right) co-ordinate changes.
The group $A$ acts on the space of germs of maps
$f: ({\mathbb C}, 0)\to ({\mathbb C}^n, 0)$ (or of curves) by
$$ (g, h)\cdot f:=g\circ f\circ h^{-1}\qquad g\in L, \quad h\in R.$$

Two germs of curves $f$ and $f'$ are equivalent if they lie in one orbit
of the $A$-action.

A germ $f$ is simple if there exists a neighbourhood of $f$ in the
space of germs which intersects only finite number of $A$-orbits.
We consider the space of germs with standard Whitney's topology:
the basis of the topology consists of coimages of open sets in a
space of $k$-jets for any $k$.
A germ $f$ is stably simple if it remaines simple after the natural
immersion ${\mathbb C}^n\hookrightarrow {\mathbb C}^N$.

A reducible curve singularity is determined by a collection of maps
$({\mathbb C}, 0)\to ({\mathbb C}^n, 0)$.

\begin{definition}
A multigerm in ${\mathbb C}^n$ is a set $F=(f_1,\ldots, f_k)$ of germs of
analytic maps $f_i: ({\mathbb C}, 0)\to ({\mathbb C}^n, 0)$, where
$\Image{f_i}\cap \Image{f_j}=\{0\}$ for $i\ne j$
($f_1$,\ldots $f_k$ are called components of the multigerm).
\end{definition}

Let $A=L\times R_{(1)}\times \ldots \times R_{(k)}$, where $R_{(i)}$ is
(the $i$-th copy of) the group of right equivalences. The group $A$
(of right-left equivalences) acts on the space of multigerms by the
formula
$$ (g, h_1,\ldots h_k)\cdot (f_1, \ldots f_k) =
     (g\circ f_1\circ h_1^{-1},\ldots, g\circ f_k\circ h_k^{-1}).$$

\begin{definition}
A multigerm $F=(f_1,\ldots, f_k)$ is called simple if there exists
a neighbourhood of $F$ in the space of multigerms which intersects
only finite number of $A$-orbits. It is stably simple, if it
remaines simple after the immersion
${\mathbb C}^n\hookrightarrow {\mathbb C}^N$.
\end{definition}

\begin{definition}
Two multigerms $F$ and $F'$ in ${\mathbb C}^n$ are equivalent if
they lie in one orbit of the $A$-action.
\end{definition}

We shall classify stably simple multigerms with respect to
the {\em stable\/} equivalence, see~\cite{Arn}.

\section{Statement of the classification}

Denote by $G_n$ a multigerm consisting of co-ordinate axes in ${\mathbb C}^n$,
more exactly
$$\begin{array}{lccr}
(\scr t_1, &  \scr 0,   & \scr \ldots &  \scr 0) \\
(\scr 0,   & \scr t_2,  & \scr \ldots &  \scr 0) \\
\scr\vdots & \scr\vdots & \scr\ddots  & \scr\vdots \\
(\scr 0,   & \scr 0,    & \scr \ldots & \scr t_n)
\end{array} $$
We denote by $(t^m\times k)$ an irreducible curve
of the form $(\underbrace{t^m,\ldots, t^m}_{k})$.

The aim of the paper is to prove the following result.

\begin{theorem}
Every stably simple multigerm up to permutation of curves
is stably equivalent to one and only
one multigerm from the following list ($m$, $n$, $k$ and $l$ are
natural numbers):
\end{theorem}

\begin{center}
\large\bf 1 Pairs of curves with regular first component
\end{center}

The first component is equal to $(t, 0)$. We write only the
second one. \\

{\bf 1.1 Multigerms with both regular components}
$$
\scriptstyle
\begin{array}{ll}
\scr 1.\; (0,t)  & \scr 2.\; (t,t^k)
\end{array}
$$

{\bf 1.2 Multigerms with the second component of multiplicity $2$}
$$
\scriptstyle
\begin{array}{lll}
\scr 1.\; (t^2,t^{2m+1}) &
\scr 2.\; (t^2,t^{2m+1}+t^{2n})\quad     m < n < 2m &
\scr 3.\; (t^2,t^{2m+1},t^{2n})\quad     m < n\le 2m \\
\scr 4.\; (t^2,t^{2m+1}+t^{2n},t^{2s})\quad m < n < s\le  2m &
\scr 5.\; (t^2,t^{2n}+t^{2m+1})\quad     n\le m &
\scr 6.\; (t^2,t^{2n},t^{2m+1})\quad     n\le m \\
\scr 7.\; (t^2,t^{2n}+t^{2m+1},t^{2s+1})\quad  n\le m < s < m+n &
\scr 8.\; (t^{2r+1},t^2) &
\scr 9.\; (0,t^2,t^{2r+1})
\end{array}
$$

{\bf 1.3 Pairs with the 3-jet $((t, 0), (0, t^3))$}
$$
\scriptstyle
\begin{array}{ll}
\scr 1.\; (t^{3m+1},t^{3})                             &
\scr 2.\; (t^{3m+1},t^{3},t^{3n+2}) \quad m\le n\le 2m \\
\scr 3.\; (t^{3m+1}+t^{3n+2}, t^3) \quad m\le n<2m    &
\scr 4.\; (t^{3m+1}+t^{3n+2},t^3,t^{3l+2}) \quad m\le n<l\le 2m \\
\scr 5.\; (t^{3m+2},t^3) &
\scr 6.\; (t^{3m+2},t^3,t^{3n+1}) \quad m<n\le 2m+1 \\
\scr 7.\; (t^{3m+2}+t^{3n+1}, t^3) \quad m<n\le 2m &
\scr 8.\; (t^{3m+2}+t^{3n+1}, t^3, t^{3l+2}) \quad m<n<l\le 2m+1 \\
\scr 9.\; (0,t^3,t^{3m+1}) &
\scr 10.\; (t^{3n+2},t^3,t^{3m+1}) \quad m\le n<2m \\
\scr 11.\; (t^{3l+2},t^3,t^{3m+1}+t^{3n+2}) \quad m\le n\le l<2m,
            \mbox{besides } n=l=2m-1 &
\scr 12.\; (0,t^3,t^{3m+2}) \\
\scr 13.\; (t^{3l+1},t^3,t^{3m+2}+t^{3n+1}) \quad m<n\le l\le 2m,
            \mbox{besides } n=l=2m &
\scr 14.\; (t^{3n+1},t^3,t^{3m+2}) \quad m<n\le 2m \\
\scr 15.\; (0,t^3,t^{3m+1},t^{3n+2}) \quad m\le n<2m &
\scr 16.\; (0,t^3,t^{3m+1}+t^{3n+2}) \quad m\le n<2m-1 \\
\scr 17.\; (0,t^3,t^{3m+1}+t^{3n+2},t^{3l+2}) \quad m\le n<l<2m &
\scr 18.\; (0,t^{3},t^{3m+2},t^{3n+1}) \quad m<n\le 2m \\
\scr 19.\; (0,t^3,t^{3m+2}+t^{3n+1}) \quad m<n<2m &
\scr 20.\; (0,t^3,t^{3m+2}+t^{3n+1},t^{3l+1}) \quad m<n<l\le 2m
\end{array}
$$

{\bf 1.4 Multigerms with the $3$-jet $((t,0),(t^3,0))$}
$$
\scriptstyle
\begin{array}{llll}
\scr 1.\;  (t^3,t^4)        &
\scr 2.\;  (t^3,t^4,t^5)    &
\scr 3.\;  (t^3,t^4,t^5,t^6)&
\scr 4.\;  (t^3,t^4+t^6)    \\
\scr 5.\;  (t^3,t^4+t^6,t^9)&
\scr 6.\;  (t^3,t^4,t^6)    &
\scr 7.\;  (t^3,t^4,t^9)    &
\scr 8.\;  (t^3,t^5,t^6)    \\
\scr 9.\;  (t^3,t^5,t^6+t^7)&
\scr 10.\;  (t^3,t^5,t^6,t^7)&
\scr 11.\;  (t^3,t^5+t^6,t^7)&
\scr 12.\;  (t^3,t^5+t^6,t^7,t^9)\\
\scr 13.\;  (t^3,t^5+t^6,t^9)&
\scr 14.\;  (t^3,t^5,t^7)&
\scr 15.\;  (t^3,t^5,t^7,t^9)&
\scr 16.\;  (t^3,t^5,t^9)\\
\scr 17.\;  (t^3,t^5+t^6,t^{12})&
\scr 18.\;  (t^3,t^5+t^6)&
\scr 19.\;  (t^3,t^5)&
\scr 20.\;  (t^3,t^5,t^{12})\\
\scr 21.\;  (t^3,t^5+t^9)&
\scr 22.\;  (t^3,t^5+t^9,t^{12})&
\scr 23.\;  (t^3,t^6,t^7,t^8)
\end{array}
$$

{\bf 1.5 Multigerms with the $4$-jet $((t,0),(0,t^4))$ or $((t,0),(t^4,0))$}
$$
\scriptstyle
\begin{array}{llll}
\scr 1.\;  (t^5,t^4,t^6,t^7)     &
\scr 2.\;  (t^6,t^4,t^5,t^7)     &
\scr 3.\;  (0,t^4,t^5,t^7)       &
\scr 4.\;  (0,t^4,t^5,t^6)       \\
\scr 5.\;  (0,t^4,t^5,t^6,t^7)   &
\scr 6.\;  (t^7,t^4,t^5,t^6)     &
\scr 7.\;  (0,t^4,t^6,t^7,t^9)   &
\scr 8.\;  (0,t^4,t^6,t^7)       \\
\scr 9.\;  (t^9,t^4,t^6,t^7)     &
\scr 10.\;  (t^4,t^5,t^6,t^7)     &
\scr 11.\;  (t^4,t^5,t^6,t^7,t^8)  
\end{array}
$$

{\bf 1.6 Curves with the $5$-jet $(0,t^5)$}
$$
\scriptstyle
\begin{array}{lllll}
\scr 1.\;  (t^9,t^5,t^6,t^7,t^8)  &
\scr 2.\;  (0,t^5,t^6,t^7,t^8)    &
\scr 3.\;  (0,t^5,t^6,t^7,t^8,t^9)&
\scr 4.\;  (t^8,t^5,t^6,t^7,t^9)  &
\scr 5.\;  (0,t^5,t^6,t^7,t^9)
\end{array}
$$
\pagebreak

\begin{center}
\large\bf 2 Pairs of curves with singular components
\end{center}

{\bf 2.1 Infinite series}\\
The first component is equal to $(t^2, t^{2m+1})$.
We write only the second component (here $m\le n$).
$$
 \begin{array}{lll}
  \scr 1.\; (t^{2n+1}, t^2) & \scr 2.\; (t^{2n+1}, t^2, t^{2n+3}) &
  \scr 3.\; (0, t^2, t^{2n+1}) \\
  \scr 4.\; (t^{2n+1}, 0, t^2) & \scr 5.\; (0, t^{2n+1}, t^2) &
  \scr 6.\; (0, 0, t^2, t^{2n+1})
 \end{array}
$$

{\bf 2.2 Individual singularities}\\
First component is equal to $(t^2, t^3)$.
$$
 \begin{array}{llll}
  \scr 1.\; (t^2, 0, t^3, t^4) & \scr 2.\; (t^2, 0, t^3) & 
  \scr 3.\; (0, 0, t^3, t^4, t^5) & \scr 4.\; (0, t^5, t^3, t^4) \\ 
  \scr 5.\; (t^5, 0, t^3, t^4) & \scr 6.\; (0, 0, t^3, t^4) &
  \scr 7.\; (0, t^4, t^3, t^5) & \scr 8.\; (t^4, 0, t^3, t^5) \\
  \scr 9.\; (0, 0, t^3, t^5, t^7) & \scr 10.\; (0, t^7, t^3, t^5) & 
  \scr 11.\; (t^7, 0, t^3, t^5) & \scr 12.\; (0, 0, t^3, t^5)
 \end{array}
$$

\begin{center}
\large\bf 3 Multigerms with regular components
\end{center}
$$
 \begin{array}{lll}
  \scr 1.\; G_n & \scr 2.\; G_n, (t\times k,0,\ldots, 0)\quad 1 < k\le n &
  \scr 3.\; G_n, (t,t^m\times k,0,\ldots, 0)\quad 1\le k<n,\; m>1 \\
  \scr 4.\; G_n, (t,0\times{(n-1)},t^m)\quad m>1 & 
  \scr 5.\; (t_1,0,0), (t_2,t_2^2,0), (t_3,0,t_3^2) &
 \end{array}
$$

\begin{center}
\large\bf 4 Multigerms with one singular component
\end{center}
 
{\bf 4.1 Series with any number of regular components}\\
The regular part is equal to $G_n,\; n\ge 2$. We write only the
singular component.
$$
 \begin{array}{ll}
  \scr 1.\; (0\times n, t^2, t^{2m+1}) & 
  \scr 2.\; (t^{2m+1}\times k, 0\times{(n-k)}, t^2) \quad 1\le k\le n \\
  \scr 3.\; (t^2\times k, 0\times{(n-k)}, t^3)\quad 1< k\le n &
  \scr 4.\; (t^2, 0\times{(n-1)},t^3,t^4) \\
  \scr 5.\; (t^2, t^4\times k, 0\times{(n-k-1)}, t^3)\quad 0\le k<n &
  \scr 6.\; (0\times n, t^3, t^4, t^5) \\
  \scr 7.\; (t^5\times k, 0\times{(n-k)}, t^3, t^4)\quad 0\le k\le n &
  \scr 8.\; (t^4\times k, 0\times{(n-k)}, t^3, t^5)\quad 0\le k\le n \\
  \scr 9.\; (0\times n,t^3, t^5, t^7) &
  \scr 10.\; (t^7\times k, 0\times{(n-k)}, t^3, t^5)\quad 1\le k\le n
 \end{array}
$$ 

{\bf 4.2 Infinite series with two regular components}\\ 
The regular part is equal to $G_2$. We write only the singular component.
$$\scriptstyle 
 \begin{array}{ll}
  \scr 1.\;  (t^2, t^2, t^{2m+1})\quad m\ge 2 & 
  \scr 2.\;  (t^2, t^2+t^{2m+1}, t^{2m+3})\quad m\ge 1 \\
  \scr 3.\;  (t^2, t^2+t^{2m+1})\quad m\ge 1 &
  \scr 4.\;  (t^2, 0, t^{2m+1}, t^{2n})\quad m<n\le 2m \\
  \scr 5.\;  (t^2, t^{2n}, t^{2m+1}+t^{2n})\quad m<n<2m &
  \scr 6.\;  (t^2, t^{2n}, t^{2m+1})\quad m<n\le 2m \\
  \scr 7.\;  (t^2, 0, t^{2m+1}+t^{2n}, t^{2l})\quad m<n<l\le 2m &
  \scr 8.\;  (t^2, t^{2l}, t^{2m+1}+t^{2n})\quad m<n<l\le 2m \\
  \scr 9.\;  (t^2, 0, t^{2m+1}+t^{2n})\quad m<n<2m &
  \scr 10.\; (t^2, 0, t^{2m+1}) \\
  \scr 11.\; (t^2, t^{2m+1}, t^{2n})\quad m<n\le 2m+1 &
  \scr 12.\; (t^2, t^{2m+1}+t^{2n}, t^{2l})\quad m<n<l\le 2m+1 \\
  \scr 13.\; (t^2, t^{2m+1}+t^{2n})\quad m<n\le 2m &
  \scr 14.\; (t^2, t^{2m+1}) \\
  \scr 15.\; (t^2, 0, t^{2m}, t^{2n+1})\quad 1<m\le n &
  \scr 16.\; (t^2, t^{2n+1}, t^{2m})\quad m\le n \\
  \scr 17.\; (t^2, 0, t^{2m}+t^{2n+1}, t^{2l+1})\quad m\le n<l<m+n &
  \scr 18.\; (t^2, t^{2l+1}, t^{2m}+t^{2n+1})\quad m\le n\le l<m+n \\
  \scr 19.\; (t^2, 0, t^{2m}+t^{2n+1})\quad 1<m\le n &
  \scr 20.\; (t^2, t^{2m}, t^{2n+1})\quad 1<m\le n \\
  \scr 21.\; (t^2, t^{2m}+t^{2n+1}, t^{2l+1})\quad m\le n<l\le m+n &
  \scr 22.\; (t^2, t^{2m}+t^{2n+1})\quad 1<m\le n
 \end{array}
$$ 

{\bf 4.3 Individual singularities} \\
The regular part is equal to $G_2$.
$$
 \begin{array}{llll}
  \scr 1.\; (t^3, t^3, t^4, t^5) & \scr 2.\; (t^3, 0, t^4, t^5, t^6) &
  \scr 3.\; (t^3, t^6, t^4, t^5) & \scr 4.\; (t^3, 0, t^4, t^5) \\ 
  \scr 5.\; (0, 0, t^4, t^5, t^6, t^7) & \scr 6.\; (t^7, t^7, t^4, t^5, t^6) &
  \scr 7.\; (t^7, 0, t^4, t^5, t^6) & \scr 8.\; (0, 0, t^4, t^5, t^6) \\
  \scr 9.\; (t^6, t^6, t^4, t^5, t^7) & \scr 10.\; (t^6, 0, t^4, t^5, t^7) & 
  \scr 11.\; (0, 0, t^4, t^5, t^7) &
 \end{array}
$$ 

{\bf 4.4 Series with the regular part $((t_1, 0), (t_2, t_2^2))$}\\
$$ 
 \begin{array}{lll}
  \scr 1.\; (0, 0, t^2, t^{2m+1}) & \scr 2.\; (0, t^{2m+1}, t^2) &
  \scr 3.\; (t^{2m+1}, 0, t^2)
 \end{array}
$$ 

\begin{center}
\large\bf 5 Multigerms with two singular components 
\end{center} 
These multigerms contain three components.
The first and the third components are $(t, 0, 0, 0)$ and
 $(0, 0, t_2^2, t_2^3)$ correspondently. In the following list $m\ge 1$.
$$ 
 \begin{array}{lll}
  \scr 1.\; (0, t_1^2, 0, 0, t_1^{2m+1}) & 
  \scr 2.\; (t_1^{2m+1}, t_1^2, 0, t_1^{2m+1})&
  \scr 3.\; (0, t_1^2, 0, t_1^{2m+1}) \\
  \scr 4.\; (t_1^{2m+1}, t_1^2, t_1^{2m+1}, 0) & 
  \scr 5.\; (0, t_1^2, t_1^{2m+1}, 0) &
  \scr 6.\; (t_1^{2m+1}, t_1^2, 0, 0)
 \end{array}
$$

We denote by $i.j.k$ the $k$-th multigerm from part $i.j$ of the list.

\section{The methods of classification}

Let $f$ be a germ of a curve in $({\mathbb C}^n, 0)$. Since $f$ is
a germ of an analytic mapping it can be represented by a
power series.

\begin{definition}
The power of the first non--zero monomial in
the power decomposition of $f$ is called the multiplicity of $f$.
\end{definition}

We shall also use the concept of invariant semigroup of $f$ 
from~\cite[section~4]{GH}.

By $M_n$ denote the ring of germs of analytic maps 
$({\mathbb C}^n, 0)\to ({\mathbb C}, 0)$. $M_n$ is ideal in
the ring of germs of analytic maps $({\mathbb C}^n, 0)\to {\mathbb C}$.

\begin{definition}
Let $f$ be a germ of a curve in $({\mathbb C}^n, 0)$. For any 
$\phi: ({\mathbb C}, 0)\to ({\mathbb C}, 0)$
there is a natural valuation $\ord: \phi\mapsto\ord{\phi}$, where 
$\ord{\phi}$ is the order of the power decomposition of $\phi$ at $0$. 
Note that $\ord$ is a semigroup homomorphism.
$S_k(f) = \ord{f^*(M^k_n)}$ is called an invariant semigroup, since
it is an $A$-invariant of $f$. $S_0(f)$ is called
the (classical) value semigroup.
\end{definition}

\begin{definition}
Let $p$ be the smallest
positive integer in the value semigroup, and
$q$ be the smallest
integer in the semigroup which is greater than $p$
and is not a multiple of $p$. We refer to $p$ as the
{\em multiplicity} of the germ, and to the pair $(p, q)$ as  the
{\em invariant} pair (see~\cite{GH}).
\end{definition}

Denote by $A_r$ the subgroup of $A$ consisting of those $A$-changes whose 
$r$-jet is equal the identity, $A_r\triangleleft A$.

We use the following statement (see~\cite[Proposition 4.2]{GH}).

\begin{lemma}
Let $f$ be a germ of a curve in $({\mathbb C}^n, 0)$ with
$n\ge 2$. Then the largest integer $N_k$ which is
not in the invariant semigroup
$S_k(f)$ is at the same time
the degree of $L_{(k-1)}$-determinacy of $f$.
\end{lemma}

Now we shall formulate a very important theorem from~\cite{BKdP}.

\begin{lemma}[Mather; see~\cite{BKdP}]
Let $G$ be a Lie group acting smoothly on a finite
dimensional manifold $V$. Let $X$ be a connected submanifold of $V$. Then
$X$ is contained in a single orbit of $G$ if and only if
\begin{enumerate}
\item for each $x\in X$ the tangent space $T_xX\subset T_x(G\cdot x)$ 
\item $\mbox{dim } T_x(G\cdot x)$ is constant for all $x\in X$.
\end{enumerate}
\end{lemma}

We shall use the statement in the following situation. $V$ is a space of
$m$-jets of multigerms with $k$ components. $G = A^{(m)}$ the group of
$m$-jets of $A$-changes. 

\begin{lemma}[The method of complete transversals]
Let $F$ be a
multigerm in $({\mathbb C}^n, 0)$
with $k$ components and let $j^mF$ be
the $m$-jet of $F$. Let $T$ be a vector subspace of the space $H^{m+1}$ of
homogeneous polynomial multigerms of degree $m+1$ with $k$ components in 
${\mathbb C}^n$ such that  
$$ T(A^{(m+1)}_1\cdot F)+T \supset H^{(m+1)}. $$
Then any $(m+1)$-jet with $m$-jet $j^mF$ is $A^{(m+1)}_1$-equivalent
to $F+t$ for some $t\in T$.
\end{lemma}

This statement can be derived from the Proposition~2.2 from~\cite{BKdP}.
The space $T$ is called a complete transversal. Sometimes the (affine)
space $F+T$ is also called a complete transversal.  

The general method is the following. We fix the 1-jet of a multigerm and
move to higher jets using the method of complete transversals. When we obtain
a complete transversal we try to simplify it using the Mather Lemma. If we
obtain finite number of jets we consider each of them separately. If our
$m$-jet is $m$-determined we add it to our list. 
If we obtain a family which can be parameterised by
a parameter (or by parameters) we conclude that this jet and all jets
adjacent to it are not simple. To prove that a given $k$-jet is not simple
we often use the following observation. We consider a special
submanifold $M\subset J^k$ that contains the jet. We prove that the tangent
space to the $A^{(k)}$-orbit does not contain
the tangent space to $M$ at any
point of it, so each point of $M$ is not simple. In particular, if
the dimension of the tangent space to $A^{(k)}$-orbit is less than
the dimension of the submanifold $M$,
then each point of $M$ is not simple.

Sometimes (as in~\cite{Arn})
we denote a germ $(t^m+t^n, t^k,\ldots)$ by $((m, n), k,\ldots)$.

\part{Multigerms with two components}

\section{Pairs of curves with one regular component}

\subsection {Multigerms with both components regular}

We refer to a curve as regular if it can be reduced to the
form $(t,0,...,0)$. In this part we assume that the first
component of a pair is regular and has
been reduced to the normal form.

\begin{lemma} \label{Lemma 1.1}
A pair with two regular components is equivalent to
one of the multigerms $1.1.1$ or $1.1.2$.
\end{lemma}

\begin{proof} $1$-jet of the second component is $(at,bt)$
which is equivalent to $(t,0)$ or $(0,t)$. Consider the first case.
Let $k$ be the minimal number such that the $k$-jet is not $(t,0)$
(if k is infinity then the pair is not simple). Then the multigerm is
$k$-determined and is equivalent to $((t,0),(t,t^k))$.
In the second case the multigerm is $1$-determined and
we obtain the normal form $1.1.1$.
\end{proof}

\subsection {Multigerms with the second component of
multiplicity $2$}

Now assume that the 1-jet of the second component is
trivial. The nontrivial $2$-jet is equvalent to $(t^2,0)$ or $(0,t^2)$.

\begin{lemma}\label{Lemma 1.2}
The second components with the $2$-jet $(t^2,0)$
or $(0,t^2)$ are equivalent to $1.2.1$--$1.2.9$.
\end{lemma}

\begin{proof} 
Suppose the $2$-jet of the second component is
$(t^2,0)$. Let $k$ be the minimal number such that the
$k$-jet is not $(t^2,0)$.

Suppose $k$ is odd and is equal to $2m+1$. Then the $k$-jet is
equivalent to
$(t^2,t^{2m+1})$. If there are no more nontrivial monomials then
we obtain the normal form $1.2.1$ since
it is at most $4m$-determined. Note that a complete
transversal in $J^n$ is trivial for odd $n$. 
Therefore a nontrivial complete transversal, for some $n\le 2m$, is
$$(t^2,t^{2m+1}+bt^{2n},ct^{2n}).$$

If $c\ne 0$ then we obtain the normal form $1.2.3$.
If, for all $k>2n$, the $k$-jet equals the $2n$-jet then
we obtain the normal form $1.2.2$, otherwise, for some $s$ in
$J^{2s}$, a complete transversal lies in
$$(t^2,t^{2m+1}+t^{2n}+bt^{2s},ct^{2s}),\ \  n<s\le 2m.$$

If $c\ne 0$ then we obtain the normal form $1.2.4$. Let
$c=0$.
The tangent space to $A^{2s}$ contains the vectors
$(2t^{2+2s-2n},2nt^{2s}+(2m+1)t^{2(m+s-n)+1})$,
$(t^{2+2s-2n},0)$, $(t^{2s},0)$ and
$(0,t^{2(m+s-n)+1}+t^{2s})$. Therefore we obtain
the normal form $1.2.2$.

Consider the case when $k$ is even and is equal to $2n$, i.e.,
when the $k$-jet is equivalent to $(t^2,t^{2n})$.
Let $2m+1$ be the minimal order of jet which is different
from $(t^2,t^{2n})$.
In $J^{2m+1}$ it is equvalent to
$$(t^2,t^{2n}+bt^{2m+1},ct^{2m+1}),\ \  m\ge n$$

If $c\ne 0$ then the second component is equivalent to 
the normal form $1.2.6$. Suppose $c=0$.
Then the $(2m+1)$-jet is equivalent to
$(t^2,t^{2n}+t^{2m+1})$. If the higher jets
are equal to the $(2m+1)$-jet then we obtain the normal form
$1.2.5$. In the other case there
exists $J^{2s+1}$ where a complete transversal is
$$(t^2,t^{2n}+t^{2m+1},ct^{2s+1}).$$

Note that the tangent space to the $A^{2m+2n+1}$-orbit
contains the vectors $(0,0,t^{4n}+2t^{2(m+n)+1})$,
$(0,0,t^{4n}+t^{2m+2n+1})$ and therefore the vector
$(0,0,t^{2m+2n+1})$.
It gives the restriction $s<m+n$.

If $c\ne 0$ then we obtain the normal form $1.2.7$.
Otherwise the second component is equivalent 
to the normal form $1.2.5$ since the tangent space to
the $A^{2s+1}$-orbit contains the vectors
$(t^{2s+1},0,0)$ and
$(2t^{2+2s-2m},2nt^{2n+2s-2m}+(2m+1)t^{2s+1})$,
$(0,t^{2n+2s-2m}+t^{2s+1})$.

Suppose that the $2$-jet is $(0,t^2)$. If, for all $k>2$, the
$k$-jet is equal to the $2$-jet then the multigerm is not simple.
Otherwise the second component is equvalent to the normal form $1.2.8$ or
$1.2.9$.
\end{proof}

\subsection {Pairs with the 3-jet $((t, 0), (0, t^3))$}

\begin{lemma}\label{Lemma 1.3}
Any simple pair of curves
with the 3-jet $((t,0),(0,t^3))$ is equivalent to one
of $1.3.1$--$1.3.20$
\end{lemma}

\begin{proof}
In $J^{3m}$ a complete transversal is trivial. In $J^{3m+1}$ it
is equivalent to $$(at^{3m+1},t^3,bt^{3m+1}). \eqno (1)$$
We can eliminate the monomials in the other co-ordinates by simple
permutation of co-ordinates and left equivalences.

If $b\ne 0$ we obtain the $(3m+1)$-jet $(0,t^3,t^{3m+1})$. This
jet is $(6m-1)$-determined. A complete transversal in $J^{3n+1}$ is
trivial. If $n<2m$ a complete transversal in $J^{3n+2}$ is
equivalent to $$(ct^{3n+2}, t^3, t^{3m+1}+dt^{3n+2}, et^{3n+2}).\eqno (2)$$
If $e\ne 0$ we obtain the normal form $1.3.15$.
If $e=0$ and $c\ne 0$ we obtain $1.3.10$ if $d=0$ or $1.3.11$
($n=l$) if $d\ne 0$. 
In the last case, if $n=l=2m-1$ we obtain
$(t^{6m-1},t^3, t^{3m+1}+t^{6m-1})$.
The tangent space to the $A^{(6m-1)}$-orbit
contains the vectors
$(0, 3t^{3m+1}, (3m+1)t^{6m-1})$ and $(0, t^{6m-1},0)$,
so, using the Mather Lemma, we obtain the normal form $1.3.10$.
If $c=0$ and $d\ne 0$ our jet is equivalent to 
$(0, t^3, t^{3m+1}+t^{3n+2})$. We have to move to higher jets. 
A complete transversal in $J^{3l+2}$, where $n<l<2m$, is equivalent to
$$(at^{3l+2}, t^3, t^{3m+1}+t^{3n+2}, bt^{3l+2})$$ since the
tangent space to $A^{(3l+2)}_1$-orbit contains the vectors
$(0,(3m+1)t^{3(m+l-n)+1}+(3n+2)t^{3l+2})$ and
$(0,t^{3(m+l-n)+1}+t^{3l+2})$. If $b\ne 0$ we obtain the
normal form $1.3.17$. If $b=0$ and $a\ne 0$ we obtain $1.3.11$. 
If $a=b=0$ for every $l$ such that $n<l<2m$
we obtain $1.3.16$. If in $(2)$ $c=d=0$, for every
$n$ such that $m\le n<2m$ we obtain $1.3.9$.
Similar reasonings with the
$(3m+2)$-jet $(0,t^3,t^{3m+2})$ produce the normal forms
$1.3.12$--$1.3.14$ and $1.3.18$--$1.3.20$.

Let's now suppose that $b=0$ and $a\ne 0$ in $(1)$. We obtain
the $(3m+1)$-jet $(t^{3m+1}, t^3)$. This jet is $(6m+2)$-determined. A
complete transversal in $J^{3n+1}$ is trivial. If $m\le n\le 2m$
then a complete transversal in $J^{3n+2}$ is equivalent to
$$(t^{3m+1}+ct^{3n+2}, t^3, dt^{3n+2}). \eqno(3)$$
If $d\ne 0$ we
obtain the normal form $1.3.2$. If $d=0$ and $c\ne 0$ we obtain
the $(3n+2)$-jet $(t^{3m+1}+t^{3n+2}, t^3)$. If $n=2m$ we have the
$(6m+2)$-jet $(t^{3m+1}+t^{6m+2}, t^3)$. Since the tangent space to
the $A^{6m+2}$-orbit contains the vectors $((t^2, 0), (t^{6m+2}, 0))$
and $((2t^2, 0), (0, 0))$ (here we write both components),
we obtain $(t^{3m+1}, t^3)$, i.e., the normal form $1.3.1$. If $n<2m$ we
move to the higher jets. Consider a complete transversal in
$J^{3l+2}$, where $n<l\le 2m$. Since the tangent space to the
$A^{3l+2}_1$-orbit contains the vectors
$((3m+1)t^{3(m+l-n)+1}+(3n+2)t^{3l+2}, 3t^{3(l-n+1)})$,
$(t^{3(m+l-n)+1}+t^{3l+2}, 0)$ and $(0, t^{3(l-n+1)})$,
a complete transversal is equivalent to
$$(t^{3m+1}+t^{3n+2}, t^3, ht^{3l+2}).$$
If $h\ne 0$ we obtain the normal form $1.3.4$. If
$h=0$, for any $l$ such that $n<l\le 2m$ we obtain $1.3.3$.
If in $(3)$ $c=d=0$, for any $n$ such that $m\le n<2m$
we obtain $1.3.1$. Similar reasonings with the $(3m+2)$-jet $(3m+2, 3)$
produce the normal forms $1.3.5$--$1.3.8$.
\end{proof}

\subsection {Multigerms with the $3$-jet $((t,0),(t^3,0))$}

We begin with the following lemma.

\begin{lemma}\label{Lemma 1.4}
A pair, the second component of which has the 8-jet $(t^3,t^6,t^7)$,
is not simple.
\end{lemma}

\begin{proof}
Let $X$ be the $11$-dimensional space
$$(a_1t^3+a_2t^4+a_3t^5+a_4t^6+a_5t^7+a_5t^8,b_1t^6+b_2t^7+b
_3t^8,c_1t^7+c_2t^8).$$
Let $x$ be a point of $X$ with nonzero co-ordinates: $a_i\ne 0,
b_i\ne 0, c_i\ne 0$. For a group $G$ acting on
the $J^8$ denote by $X_G$ the intersection $T_xX\cap T_x(Gx)$,
where $T_xX$ is the tangent space to $X$ and $T_x(Gx)$ is the
tangent space to the $G$-orbit of the point $x$. Let's evaluate the
dimension of $X_{A^8}$. The group $A^8$ is the product of
$L^8$ and $R^8$. If $l\circ x\notin X$, where $l$ is a $L^8$-change
and $x$ is a point of $X$ then, for any $R^8$-change $r$, the
point $l\circ x\circ r^{-1}$ lies out of X. Moreover if $l\circ
x\circ r^{-1}$ is a point of $X$ then the points $x\circ r^{-1}$
and $l\circ x$ lie in $X$ too. It means that
a vector $v$ from $X_{A^8}$ is the sum of some vectors
$v_1\in X_{R^8}$ and $v_2\in X_{L^8}$.
Therefore $X_{A^8}$ is the sum of $X_{R^8}$
and $X_{L^8}$. $X_{R^8}$ is the $6$-dimensional space and
includes the vectors $(t^6,0,0), (t^7,0,0),$ and $(t^8,0,0)$, since these
vectors are the linear combinations of
$(3a_1t^6+4a_2t^7+5a_3t^8,0,0), (3a_1t^7+4a_2t^8,0,0)$ and
$(3a_1t^8,0,0)$. The preimage of $X_{L^8}$ with respect to
the tangent mapping $g:T_eL^8\to T_xJ^8$
is $$<(x,0,0),(y,0,0),(z,0,0),(0,y,0),(0,0,z),(0,z,0),(x^2,0,0)>.$$
Therefore the intersection $X_{R^8}\cap X_{L^8}$ is at least
$3$-dimensional (it contains the vectors $g(x^2,0,0)$, $g(y,0,0)$,
$g(z,0,0)$) and the dimension of $X_{A^8}$ is not more than $10$.
This implies that any open neighbourhood of a point $x$
intersects with infinite number of different orbits.
\end{proof}

\begin{lemma}\label{Lemma 1.5}
Second components of simple multigerms with the
$4$-jet $(t^3,t^4)$ are classified by $1.4.1$--$1.4.7$.
\end{lemma}

\begin{proof}
First note that the second component
is at most $9$-determined and a complete transversal in
$J^7$ and $J^8$ is trivial. A complete transversal in $J^5$
is $(t^3,t^4+bt^5,ct^5)$. If $c\ne 0$ then the $5$-jet is 
equivalent to $(t^3,t^4,t^5)$. Therefore we obtain the normal forms
$1.4.2$ and $1.4.3$.
Suppose $c=0$. Then the $5$-jet is equivalent to
$(t^3,t^4)$. In $J^6$ a complete transversal is $(t^3,t^4+bt^6,ct^6)$. 
If $c\ne 0$ then we obtain the normal form $1.4.6$. Suppose $c=0$ and
$b\ne 0$ then the second component is equivalent to the normal forms
$1.4.4$ or
$1.4.5$. If $c=0$ and $b=0$ then there are two
possibilities: the normal form $1.4.1$ or $1.4.7$.
\end{proof}

\begin{lemma}\label{Lemma 1.6}
Second components of simple multigerms with the
$5$-jet $(t^3,t^5)$ are classified by $1.4.8$--$1.4.22$.
\end{lemma}

\begin{proof}
A complete transversal in $J^6$ is
$$(t^3,t^5+bt^6,ct^6).$$

Consider $c\ne 0$. Then the multigerm is 7-determined
and we obtain the normal forms $1.4.8$--$1.4.10$
Suppose $c=0$ and $b\ne 0$. A complete transversal in $J^7$
is $(t^3,t^5+t^6+et^7,ft^7)$ which is equivalent to
$(t^3,t^5+t^6,ft^7)$ since
the tangent space to the $A^7$-orbit contains the vectors
$(t^7,0,0)$, $(3t^5,5t^7,0)$,
$(t^5+t^6,0,0)$ and $(t^6,0,0)$. If $f\ne 0$ then the
multigerm is 9-determined and
we obtain the normal forms $1.4.11$ and $1.4.12$.
Suppose $f=0$. A complete transversal is not
trivial only in $J^9$ and $J^{12}$.
Since $(t^3,t^5+t^6+\alpha t^9)$ is equivalent to $(t^3,t^5+t^6)$
(the tangent space to $A^9$ includes the vectors
$(0, 5t^8+6t^9)$ and $(0, t^8)$)
we obtain the normal forms $1.4.13$, $1.4.17$ and $1.4.18$.

Consider the  case $c=0$ and $b=0$. As above a complete
transversal in $J^7$ is equivalent to $(t^3,t^5+ht^7,gt^7)$.
If $g\ne 0$ then we obtain the forms $1.4.14$ and $1.4.15$.
Suppose $g=0$, then
the $7$-jet is equvalent to $(t^3,t^5)$. Note that a complete
transversal is
not trivial only in $J^9$ and $J^{12}$. In $J^9$ it is
equivalent
to $(t^3,t^5+ht^9,lt^9)$. If $l\ne 0$ we obtain the
normal form $1.4.16$. Consider $l=0$. There are two cases: $h=0$
and $h\ne 0$.
In the first case a complete transversal in $J^{12}$ is equivalent
to $(t^3,t^5,lt^{12})$  and we obtain two normal forms:
$1.4.19$ and $1.4.20$. In the second case we obtain the $9$-jet
$(t^3,t^5+t^9)$ which produces the normal forms $1.4.21$ and $1.4.22$.
\end{proof}

Note that from Lemma~\ref{Lemma 1.4} it
follows that if the $6$-jet of the second component is
equal to $(t^3,t^6)$ then the
component is equivalent to $(t^3,t^6,t^7,t^8)$. Therefore one has

\begin{lemma}\label{Lemma 1.7}
The only multigerm with the
second component's
$6$-jet $(t^3,t^6)$ is equivalent to
$((t,0),(t^3,t^6,t^7,t^8))$, i.e., to $1.4.23$.
\end{lemma}

\subsection {Multigerms with the $4$-jet $((t,0),(t^4,0))$ or
$((t,0),(0,t^4))$}

The list is limited by the lemma.

\begin{lemma}\label{Lemma 1.8} Multigerms with the 7-jet
$$(a_1t^4+...+a_4t^7,b_1t^4+...+b_4t^7,c_1t^4+...+c_4t^7).$$
of the second component are not simple.
\end{lemma}

\begin{proof}
As in the Lemma~\ref{Lemma 1.4} consider the
$12$-dimensional
space $X$:
$$(a_1t^4+...+a_4t^7,b_1t^4+...+b_4t^7,c_1t^4+...+c_4t^7).$$
In a point $x$ with nonzero co-ordinates $a_i\ne 0, b_i\ne 0, c_i\ne 0$
the dimension of $X_{L^7}$ is $7$ and the dimension of $X_{R^7}$ is $4$.
Therefore the dimension of $X_{A^8}$ is not more than $11$. From
this it follows that, for any neibourhood $U(x)$ of $x$,
the number of orbits which intersects $U(x)$ is infinite.
\end{proof}

\begin{lemma}\label{Lemma 1.9}
The second components of simple multigerm with
$4$-jet $(0,t^4)$ are classified by $1.5.1$--$1.5.9$.
\end{lemma}

\begin{proof}
In $J^5$ a complete transversal is $$(at^5,t^4,bt^5).$$
Let $b\ne 0$. A complete transversal in $J^6$ is
$(ct^6,t^4,t^5+dt^6,ft^6).$
If $f=0$, then from Lemma~\ref{Lemma 1.8} it
follows that there are two simple components
$(t^6,t^4,t^5,t^7)$ and $(0,t^4,t^5,t^7)$. Consider $f\ne 0$.
Then the $6$-jet is equivalent to $(0,t^4,t^5,t^6)$. It produces the forms
$1.5.4$, $1.5.5$ and $1.5.6$.

Let $a\ne 0$ and $b=0$. Then the $5$-jet is equivalent to $(t^5,t^4)$ and
a complete transversal in $J^6$ lies in
$$(t^5+dt^6,t^4,ft^6).$$ Note that if $f=0$, then the pair
is adjoint to the pair of Lemma~\ref{Lemma 1.8}
and therefore is not simple. So the $6$-jet is equvalent to
$(t^5,t^4,t^6)$. From Lemma~\ref{Lemma 1.8} we obtain that
a simple component with such a $6$-jet is equivalent to $(t^5,t^4,t^6,t^7)$.

Let $a=0, b=0$. A complete transversal in $J^6$ lies in
$(ct^6,t^4,dt^6)$. Suppose
$d=0$, then the component
is adjoint to the component of the Lemma~\ref{Lemma 1.8}. So
$d\ne 0$ and the $6$-jet is equvalent to $(0,t^4,t^6)_6$.
From Lemma~\ref{Lemma 1.8} it follows that the $7$-jet
is equivalent to $(0,t^4,t^6,t^7)$, what produces the last three forms.
\end{proof}

Now consider multigerms with the second component $(t^4,0)_4$.
From Lemma~\ref{Lemma 1.8} it follows that its
$7$-jet is equivalent to $(t^4,t^5,t^6,t^7)$. It is
$8$-determined, therefore we obtain two multigerms more:

\begin{lemma}
The second components of a simple multigerm
with $4$-jet $(t^4,0)$ has one of two forms $1.5.10$ or
$1.5.11$.
\end{lemma}

\subsection {Curves with the $5$-jet $(0,t^5)$ or $(t^5,0)$}

First consider the multigerms with the 5-jet $((t, 0),(0,t^5))$.

\begin{lemma}\label{Lemma 1.10}
Pair of curves with the $9$-jet $((t,0), (t^8,t^5,t^6,t^7))$ is not simple.
\end{lemma}

\begin{proof}
The dimension of the space $X$:

$$(a_1t^8+a_2t^9,b_1t^5+b_2t^6+\ldots +b_5t^9,c_1t^6+\ldots +c_4t^9,
d_1t^7+\ldots +d_3t^9)$$ is $14$. At the point with nonzero co-ordinates
the dimensions of $X_{R^9}$ is not more than $5$ and that of
$X_{L^9}$ is $7$. Therefore
the dimension of $X_{A^9}$ is less than $14$ and the pair
is not simple.
\end{proof}

\begin{lemma}\label{Lemma 1.11}
The second components of a simple multigerm with the
$5$-jet $(0,t^5)$ has one of the forms $1.6.1$--$1.6.5$.
\end{lemma}

\begin{proof}
From Lemma~\ref{Lemma 1.8} it follows that the $7$-jet of
the second component is equivalent to $(0,t^5,t^6,t^7)$.
A complete transversal in $J^8$ lies in
$$(at^8,t^5,t^6+bt^8,t^7+ct^8,dt^8).$$
If $d\ne 0$ then we obtain the first three normal forms.
Otherwise the $8$-jet is equivalent to $(at^8,t^5,t^6,t^7)$.
Lemma~\ref{Lemma 1.10} implies that the $9$-jet is equivalent to
$(at^8,t^5,t^6,t^7,t^9)$ and we obtain the normal forms
$1.6.4$ and $1.6.5$.
\end{proof}

Multigerms with the $5$-jet $((t,0),(t^5,0))$ are not simple
since they are adjoint to the multigerm of Lemma~\ref{Lemma 1.8}.
By the same reason there are no more simple pairs with
regular first component.

\section{Pairs of curves with singular components}

Here we denote parameters on both components by $t$.

\begin{lemma}\label{Lem2.1}
Pair of curves with the $3$-jet $((t^2, t^3), (t^2, \alpha t^3))$,
where $\alpha \ne 1$, is not simple.
\end{lemma}

\begin{proof}
Let 's consider the tangent space to the $A^{(3)}$-orbit. It is generated 
by the following $8$ vectors: $((t^2, 0), (t^2, 0))$, $((0, t^2), (0, t^2))$,
$((t^3, 0), (\alpha t^3, 0))$, $((0, t^3), (0, \alpha t^3))$, $((2t^2, 3t^3),
(0,0))$, $((t^3, 0), (0, 0))$, $((0, 0), (2t^2, 3\alpha t^3))$, 
$((0, 0), (2t^3, 0))$. Note, that $$ ((0, 0), (2t^2, 3\alpha t^3))-
2((t^2, 0), (t^2, 0))+((2t^2, 3t^3), (0, 0))=3((0, t^3), (0, \alpha t^3)).$$  
So, these vectors are linearly dependent and we can remove 
$((0, t^3), (0, \alpha t^3))$ from the list above. Hence it is obvious that  
the tangent space to the $A^{(3)}$-orbit does not
contain the vector $((0, 0), (t^3, 0))$ from the tangent space to our 
1-di\-men\-sional submanifold. But if the tangent space to the
$A^{(3)}$-orbit does not contain the tangent space to the submanifold, 
for each point of it, the multigerm fails to be simple.
\end{proof}

Let us consider the 2-jet $((t^2, 0), (0, t^2))$ first. 

\begin{lemma}
Any simple multigerm with the 2-jet of the second component $(0, t^2)$, such that
the first component equals $(2, 2m+1)$ and the invariant
pair of the second component is greater than or equals $(2, 2m+1)$,
is equivalent to one from the forms $2.1.1$--$2.1.3$.
\end{lemma}

\begin{proof}
We shall not calculate tangent spaces in obvious cases, but in other cases we have
to do it. In $J^{2m}$ a complete transversal is trivial for any $m$. In $J^{2m+1}$
it is equivalent to $$((t^2, at^{2m+1}, 0), (bt^{2m+1}, t^2, ct^{2m+1})).$$ 
There exists $m$ such that $a\ne 0$ (we suppose that the multiplicity of the first
component is less than or equals the second one).  If $c\ne 0$
then we obtain $((t^2, t^{2m+1}), (0, t^2, t^{2m+1}))$. It is clear,
that this multigerm is $(2m+1)$-determined, i.e., we obtain $2.1.3$, $n=m$. 
Now suppose $c=0$.
If $b\ne 0$ then using multiplications of 
co-ordinates in the target and parametrs in the sources by constants we obtain
$$((t^2, t^{2m+1}), (t^{2m+1}, t^2)).\eqno (4)$$ If $b=0$ we 
obtain $$((t^2, t^{2m+1}), (0, t^2)).\eqno (5).$$ We shall discuss (4) later.
Now consider (5). In $J^{2n}$ a complete transversal is trivial for any $n$.
In $J^{2n+1}$ a complete transversal is given by
$$((t^2, t^{2m+1}, 0), (bt^{2n+1}, t^2, ct^{2n+1})).$$ 
If $c\ne 0$ then
we obtain $2.1.3$. If $c=0$ and $b\ne 0$ we obtain the $(2n+1)$-jet
$$((t^2, t^{2m+1}), (t^{2n+1}, t^2)).$$
Starting from this moment, it is not important for
our consideration that $n>m$, so we can also consider (4).
A complete transversal in $J^{2n+2}$ is trivial. Let 's consider it
in $J^{2n+3}$. We can simply obtain the following vectors:
$((2t^{2n+3}, 0), (0, 0))$ and $((0, 0), (0, 2t^{2n+3}))$.
Note, that the tangent space contains the vector 
$((0, t^{2n+3}), ((0, t^{2+(2n+1)(n-m+1)}))$, the power of the
last monomial is greater than or equals $2n+3$. So, the tangent
space contains $((0, t^{2n+3}), (0, 0))$. Now, we want to obtain
the vector $((0, 0), (t^{2n+3}, 0))$. It is obvious, if we note that
the tangent space contains $((t^{2m+3}, 0), (t^{2n+3}, 0))$,
$((2t^{2m+3}, (2m+1)t^{4m+2}), (0, 0))$ and 
$((0, t^{4m+2}), (0, t^{(2m+1)(2n+1)}))$,
the power of the last monomial is greater than $2n+3$. After it
we can conclude that a complete transversal is equivalent to
$$((t^2, t^{2m+1}), (t^{2n+1}, t^2, at^{2n+3})). \eqno (6)$$   
Now we shall prove, for every $a,$ that this $(2n+3)$-jet is
$L$~-determined. If, using left equivalence, we eliminate a monomial 
with even degree in one component, then we obtain a monomial of higher
degree in another component, so we have to eliminate monomials
of odd degree.

Let 's consider the $(2n+2k+1)$-jet $(k\ge 2)$ of our multigerm. 
We shall eliminate the
monomial $t^{2n+2k+1}$ in the second component. Our construction
does not depend on $a$ in (6). We can obtain this monomial as $x_1 x_2^k$.
In the first component we obtain $t^{(2m+1)k+2}$. Now we have two
possibilities: $k=2i$ or $k=2i+1$. 

In the first case we consider our multigerm in $J^{2n+4i+1}, i\ge 1$ and our 
monomial in the first component is $t^{2i(2m+1)+2}$. It can be
obtained as $x_1^{(2m+1)i+1}$. So, in the second component we obtain
a monomial of degree $(2n+1)((2m+1)i+1)=2n+1+(2n+2m+1+2nm)i>2n+1+4i$.  

In the second case we consider our multigerm in $J^{2n+4i+3}, i\ge 1$. 
The monomial in the first component is $t^{2i(2m+1)+2+2m+1}$, it can be
obtained as $x_1^{(2m+1)i+1}x_2$. So, in second component we obtain
a monomial of degree $(2n+1)((2m+1)i+1)+2=(2n+2m+2mn+1)i+2n+3>2n+4i+3$.

In this construction we did not use that $n\ge m$, so it is
suitable for eleminating the monomials in the first component as well. 

Therefore, in (6), if $a\ne 0$ we obtain the normal form $2.1.2$, if $a=0$
we obtain the normal form $2.1.1$.
\end{proof}

Now we suppose that the first component is $(t^2, t^{2m+1})$ and
the invariant pair of second component is greater than or equals $(2, 2m+1)$.

\begin{lemma}
Any pair of curves with the $(2m+1)$-jet $((2, 2m+1), (0, 0, 2))$ is
equivalent to one of the forms $2.1.4$--$2.1.6$.
\end{lemma}

\begin{proof}
In $J^{2n+1}$ a complete transversal is equivalent to
$$ ((t^2, t^{2m+1}, 0, 0), (at^{2n+1}, bt^{2n+1}, t^2, ct^{2n+1})).$$
If $c\ne 0$ then we obtain the normal form $2.1.6$ which is $(2n+1)$-determined.
If $c=0$ and $b\ne 0$ we have 
$$((t^2, t^{2m+1}, 0), (at^{2n+1}, t^{2n+1}, t^2)).$$
One can obtain the vector $((0,0,0),(t^{2n+1},0,0))$ as a linear combination of
$((t^{2m+1},0,0), (t^{2n+1},0,0))$, $((2t^{2m+1},(2m+1)t^{4m},0), (0,0,0))$ and
$((0,t^{4m},0), (0,a^{2m}t^{2m(2n+1)},0))$. Then we can use the Mather Lemma
and obtain the normal form $2.1.5$, which is $(2n+1)$-determined. If $b=c=0,\: 
a\ne 0 $ we obtain the normal form $2.1.4$ which is also $(2n+1)$-determined. 
\end{proof}

Now we shall consider the 2-jet $((2,0), (2,0))$. 

\begin{lemma}
There exist two simple curves with the 2-jet $((2,0), (2,0))$: $2.2.1$ 
and $2.2.2$.
\end{lemma}

\begin{proof}
Using Lemma~\ref{Lem2.1}, we have only one possibility for the 3-jet
$((t^2,t^3,0), (t^2,0,t^3))$. Any 4-jet of this curve is equivalent to
$$((t^2, t^3, 0, 0), (t^2+at^4, bt^4, t^3+ct^4, dt^4).$$
The tangent space to $A^{(4)}$-orbit contains the vectors $((0,0),(2t^4,0))$ 
(so we can suppose $a=0$), $((0, t^4), (0, t^4))$, $((2t^3, 3t^4),(0,0))$,
$((t^3,0),(bt^4,0))$, $((0,0), (2t^4,0))$ (so we can suppose $b=0$),
$((0,0,0), (2t^3,0,3t^4))$, $((0,0,0), (t^3+ct^4,0,0))$ and 
$((0,0,0), (2t^4,0,0))$ (so we can suppose $c=0$). Hence, if $d\ne 0$
we obtain the first curve, if $d=0$ we obtain the second curve. Note, that
both this 4-jets are 4-determined.
\end{proof}

\begin{lemma}\label{Lem2.5}
The family of 5-jets
$$((t^2, t^3, 0), (t^5, t^4+\alpha t^5, t^3))$$ is not simple.
\end{lemma}

\begin{proof}
Simple calculations in $J^5$ show 
that the tangent space to the $A^{(5)}$-orbit does not contain the vector
$((0,0,0),(0,t^5,0))$, so $\alpha$ is a module in our family. 
\end{proof}

Now we have to consider the case when the multiplicity of the second curve
equals 3 and the first curve is $(t^2, t^3)$. Curves from the other cases
are adjacent to the family from Lemma~\ref{Lem2.1}. More exactly, we have to
consider the 3-jet $((t^2,t^3,0),(0,0,t^3))$.

\begin{lemma}
Every simple pair with the 3-jet as above is equivalent to
one of $2.2.3$--$2.2.12$.
\end{lemma}

\begin{proof}
A complete transversal in $J^4$ is equivalent to
$$((t^2, t^3), (at^4, bt^4, t^3, ct^4))$$
Let 's suppose that $c\ne 0$. Then our 4-jet is equivalent to 
$((t^2,t^3,0,0),(0,0,t^3,t^4))$.
This jet is 5-determined. The tangent space to its $A^{(5)}$-orbit contains
the vectors $((0,0,0,0), (0,0,t^5,0))$, $((0,0,0,0), (0,0,3t^4,4t^5))$ and
$((0,0,0,0), (0,0,t^4,0))$, so we obtain the curves $2.2.3$--$2.2.6$.

Now let us consider the case $c=0$. If $b\ne 0$ this 4-jet is equivalent to
$$((t^2, t^3, 0), (0, t^4, t^3)). \eqno(7)$$ 
If $b=0,\: a\ne 0$ this 4-jet is equivalent to 
$$((t^2, t^3, 0), (t^4, 0, t^3)). \eqno(8)$$
If $a=b=0$ we obtain 
$$((t^2, t^3, 0), (0, 0, t^3)). \eqno(9)$$
Now we have to consider the 5-jet. Using Lemma~\ref{Lem2.5} we see that
in the second component the fourth co-ordinate equals $t^5$.
In (7) and (8) we obtain the normal forms $2.2.7$ and $2.2.8$,
which are 5-determined. 
In (9) we see that a complete transversal in $J^6$
is trivial and this curve is 7-determined. So we obtain the normal forms 
$2.2.12$, $2.2.11$, $2.2.10$ and $2.2.9$.
\end{proof}

\part{Multigerms with three and more components}

\begin{lemma}\label{Lem3_0.1}
Simple multigerm can not contain three nonregular components.
\end{lemma}

\begin{proof}
 A non-regular curve is ajacent to the curve $(2,3)$.
 That is why the first curve is adjacent to $(2, 3)$, the second one to
 $(0,0,2,3)$ and the third one to $(0,0,0,0,2,3)$.
 Note that the third curve in the triple is ajacent to $(0,0,1,0,2,3)$
 therefore
 the triple is ajacent to $((2,3,0,\ldots),(0,0,2,3,4,0),(0,0,1,0,0,3)).$
 By changing indexes of axis and curves we can obtain the triple
 $((1,0,0,0,0,3),(0,2,3,0,0,0),(2,0,0,3,4,0)).$ The second and the third
 curves are ajacent to curves with the 2-jet $(1,2)$.  So it is remaines to
 prove that a triple with the $2$-jet
 $((1,0),(1,2),(1,2))$ is not simple. Therfore we have to prove the following
 statement.
\end{proof}

\begin{lemma}\label{Lem3_0.2}
A triple with the 2-jet $((1,0),(1,2),(1,2))$ is not simple.
\end{lemma}

\begin{proof}
  Consider the $4$-dimensonal subspace of the space of 2-jets:
  $$((t,0),(\alpha_1t_1,\alpha_2t_1^2),(\alpha_3t_2,\alpha_4t_2^2)), 
                                                         \alpha_i\ne 0$$
  Using $R^{(2)}$-changes each such triple can be reduced to the same form
  with $\alpha_1= \alpha_2=1.$ Note, that if one can reduce the 2-jet to the
  form, where $\alpha_3=\beta_1, \alpha_4=\beta_2$, then one can do it 
  using only the following changes: $\tilde x=ax$, $\tilde y=by$, 
  $t_1=c\tilde t_1$ and $t_2=d\tilde t_2$. Hence we have the following 
  equations for $a, b, c ,d$:
  $$ac = ad = 1,\quad \alpha_3 bc^2=\beta_1,\quad \alpha_4 bd^2 = \beta_2.$$
  Therefore we have $c=d=a^{-1}$. Then $\alpha_3 ba^{-2}=\beta_1$ and
  $\alpha_4 ba^{-2}=\beta_2$.
  And it implies that $\alpha_3 /\beta_1 = a^2 b^{-1} = \alpha_4 /\beta_2 $. 
  So the 2-jet under consideration has a continious invariant $\alpha_3 /\alpha_4$, 
  i.e., it is not simple.
\end{proof}

\section{Multigerms with regular components}

\begin{lemma}\label{Lem3_1.1}
Any simple multigerm with all regular components is equivalent to
one of $3.1$--$3.5$.
\end{lemma}

\begin{proof}
Consider the 1-jet of our multigerm. 
Suppose that first $n$ curves can be (and are) reduced to the form $G_n$
and the 1-jets of all other curves have zero co-ordinates with number
greater than $n$ (otherwise one can reduce $n+1$ curves to the form
$G_{n+1}$). If the total number of curves is equal to $n$, we obtain the normal form 3.1.

Suppose that $n\ge 2$. Then our multigerm can not contain more than
$(n+1)$ components. The reason is that $(n+2)$-lines in ${\mathbb C}^n$
is not a simple multigerm. For $n=2$ they have a continious 
invarian~--- double ratio. For $n>2$ consider the space of 1-jets of $k$ lines 
in ${\mathbb C}^n$. It is $nk$-dimensional. The dimension of the tangent 
space to the $L^{(1)}$-orbit is less than or equals $n^2$. The dimension of the 
tangent space to the $R^{(1)}$-orbit is less than or equals $k$. 
For our multigerm to be simple it is necessary to have $nk\le n^2+k$ i.e., 
$k\le n+1+\frac{1}{n-1}$. Hence $k\le n+1$. 

So, we have $G_n$ and a non-singular curve with zero co-ordinates
with numbers greater than $n$. By $R^{(1)}$-transformations and permutations of 
curves and co-ordinates we can reduce the 1-jet to $G_n$ and \\
$(\underbrace{t,\ldots, t}_{k},0,\ldots, 0)$, where $1\le k\le n$.
If $k\ge 2$ we obtain the normal form $3.2$ (since $x_{1}x_{2}^{m}$ equals
$t_{n+1}^{m+1}$ for the last curve and zero for first $n$ curves, $m\ge 1$).
If $k=1$ we move to higher jets. A complete transversal in $J^m$ is equivalent to 
$$G_n, (t, a_1t^m,\ldots, a_nt^m).$$
If $a_n\ne 0$ we obtain the normal form $3.4$. If $a_n=0$ and there exists 
$j$ such that $a_j\ne 0$ we obtain the normal form $3.3$ (since it is $m$-determined).

Now suppose $n=1$. There are at least 3 components in our multigerm.
If the multigerm has at least 4 components, then it fails to be simple, since
its 1-jet is adjacent to 4 lines in ${\mathbb C}^2$. So, the multigerm
contains 3 components. Note, that by lemma~\ref{Lem3_0.2}
the family $((t_1, 0), (t_2, t_2^2), (t_3, \alpha t_3^2))$ is not simple. So we have 
only one possibility for the 2-jet:
$((t_1, 0, 0),(t_2, (t_2^2, 0),(t_3, 0, t_3^2))$.This jet is 2-determined and 
we obtain the normal form $3.5$.
\end{proof}

\section{Multigerms with one non--regular component}

We can suppose that the regular part of the multigerm is
equivalent to one of 5 forms from the previous section (see Lemma~\ref{Lem3_1.1})
or to $((t_1, 0), (t_2, t_2^m))$. The forms $3.2$--$3.4$ are not suitable, 
since if we add to it one non--regular curve, the 1-jet of the multigerm
will be adjacent to $(n+2)$ lines in ${\mathbb C}^n$, which is not simple.
If we add one non-regular curve to the normal form $3.5$,the 1-jet of the 
multigerm would adjoin to 4 lines in ${\mathbb C}^2$. So we have 2
possibilities for the regular part: $G_n$ or $((t_1, 0), (t_2, t_2^m))$.
First consider the case, when the multiplicity of the singular
component equals 2 and the regular part of the multigerm is $G_n$.

\begin{lemma}
Any simple multigerm with the 2-jet
$(G_n, (\underbrace{0,\ldots, 0}_{n}, t^2))$ is equivalent to 
$4.1.1$ or $4.1.2$.
\end{lemma}

\begin{proof}
A complete transversal in $J^{2m}$ is trivial for any $m$. 
A complete transversal in $J^{2m+1}$ is equivalent to
$$G_n,\quad (a_{1}t^{2m+1},\ldots, a_{n}t^{2m+1}, t^2, a_{n+1}t^{2m+1}).$$
If $a_{n+1}\ne 0$ we obtain the normal form $4.1.1$. If $a_{n+1}=0$ but there
exists $j$ such that $a_j\ne 0$ we obtain the normal form $4.1.2$ ($k$ depends on
the number of coefficients $a_j$ different from zero).
\end{proof}

Now, we shall consider other 2-jets: 
$G_n, (\underbrace{t^2,\ldots, t^2}_{k}) \quad 1\le k\le n$.
This 2-jet is not enough, so we move to 3-jets.
A complete transversal in $J^3$ is equivalent to
$$(G_n, (t^2+a_{1}t^3,\ldots, t^2+a_{k}t^3,a_{k+1}t^3,\ldots, a_{n+1}t^3)).
        \eqno (10)$$ 
First, consider $a_{n+1}\ne 0$.

\begin{lemma}\label{Lem3_2.2}
Any simple multigerm with the 3-jet
$G_n, (t^2,\ldots, t^2,\underbrace{0,\ldots, 0}_{n-k},t^3)$ is equivalent to
one of $4.1.3$--$4.1.5$. 
\end{lemma}

\begin{proof}
The normal form $4.1.3$ is 3-determined since $x_1 x_2=t^4$.
Consider the 3-jet $(G_n, (t^2,\underbrace{0,\ldots, 0}_{n-1},t^3)$. It is
5-determined so we have to consider the 4-jet. A complete transversal in $J^4$ 
is equivalent to
$$(G_n, (t^2, b_{2}t^4,\ldots, b_{n}t^4, t^3+b_{n+1}t^4, b_{n+2}t^4)).$$
If $b_{n+2}\ne 0$ we obtain the normal form $4.1.4$ which is 4-determined.
If $b_{n+2}=0$ we can eliminate $t^4$ in the first and $(n+1)$-st co-ordinates,
since the tangent space to the $A^{(4)}$-orbit contains the vectors:
$(2t^4,0,\ldots, 0)$, $(2t^3,\underbrace{0\ldots, 0}_{n-1},3t^4)$ and
$(t^3,0,\ldots, 0)$. Hence we obtain the normal form $4.1.5$ ($k$ depends 
on the number of coefficients $b_j$ different from zero).
\end{proof}

Now, suppose that $a_{n+1}=0$ in $(10)$. We shall proof that $n$ can nott be 
greater than 2. 

\begin{lemma}\label{Lem3_2.3}
The 3-jet $G_n, (a_{1}t^2+b_{1}t^3,\ldots, a_{n}t^2+b_{n}t^3)$ is not simple if $n>2$.
\end{lemma}

\begin{proof}
The dimension of the submanifold of such jets (in $J^{(3)}$)equals $2n$. 
Now we shall estimate the dimension of the stabilisator of the regular part i.e., 
the subgroup in $A^{(3)}$ which preserves the regular part. The dimension of the 
tangent space to the $L^{(3)}$-orbit is less than or equals $n$, since it is 
generated by the images of the vectors
$(\underbrace{0,\ldots, 0}_{j-1},x_j,0,\ldots, 0)$ with $1\le j\le n$.
The dimension of the tangent space to the $R^{(3)}$-orbit is less than or equals $2$,
since it is generated by the images of two vectors: $t$ and $t^2$. Hence the dimension 
of the tangent space to the $A^{(3)}$-orbit is less than or equals $(n+2)$.
If the dimension of the tangent space is less than the dimension of
the submanifold then this 3-jet fails to be simple. Hence $n+2\ge 2n$ i.e., $n\le 2$.
\end{proof}

So, if in $(10)$ $a_{n+1}=0$ then $n=2$. Consider this case in detail.
First suppose that $k=n=2$ in $(10)$.

\begin{lemma}
Any simlpe multigerm with the 3-jet $(10)$, where $k=n=2$, is
equivalent to one of $4.2.1$--$4.2.3$. 
\end{lemma}

\begin{proof}
We shall describe all simple multigerms with the 2-jet
$(G_2, (t^2, t^2))$. A complete transversal in $J^{2m}\quad (m>1)$ is 
trivial for any $m$. In $J^{2m+1}$ it is equivalent to 
$(G_2, (t^2, t^2+at^{2m+1}, bt^{2m+1}))$
since the tangent space to the $A^{(2m+1)}_1$-orbit contains the vector 
$(2t^{2m+1}, 2t^{2m+1})$. If $b\ne 0$ we obtain the normal
form $4.2.1$ since it is $(2m+1)$-determined. If $b=0$ and $a\ne 0$
we obtain the $(2m+1)$-jet $(G_2, (t^2, t^2+t^{2m+1}))$ which is 
$(2m+3)$-determined, since $x_1x_2^2-x_1^2x_2$ is equal to 
$t^{2m+5}+higher~order~terms$ for
the singular component and to zero for both regular components.
Hence we have to move to the $(2m+3)$-jet. A complete transversal in 
$J^{(2m+3)}$ is equivalent to 
$$(G_2, (t^2, t^2+t^{2m+1}, ct^{2m+3})).$$
This is so since the tangent space to the $A^{(2m+3)}_1$-orbit contains the 
vectors $v_1=(2t^{2m+3}, 2t^{2m+3})$, $v_2=(2t^4, 2t^4+(2m+1)t^{2m+3})$,
$v_3=(t^4+t^{2m+3}, 0)$ and $v_4=(0, t^4+t^{2m+3})$.
$v_1+v_2-2v_3-2v_4$ is equal to $(0, (2m+1)t^{2m+3})$ for the singular component
and to zero for both regular components.
If $c\ne 0$ we obtain the normal form $4.2.2$, if $c=0$ we obtain
the normal form $4.2.3$. 
\end{proof}

{\bf Remark:} If $m=1$ in the normal form $4.2.1$ we obtain the normal 
form $4.1.3$ with $k=n=2$.

Now we shall consider $n=2$, $k=1$ in $(10)$.

\begin{lemma}
Any simple multigerm with the 2-jet $(G_2, (t^2, 0))$ is
equivalent to one of $4.2.4$--$4.2.22$.
\end{lemma}

\begin{proof}
Consider a complete transversal in $J^{2m+1}$.
It is equivalent to 
$$(G_2, ((t^2, at^{2m+1}, bt^{2m+1})).\eqno (11)$$ 
Suppose $b\ne 0$. We obtain the $(2m+1)$-jet $(t^2, 0, t^{2m+1})$. 
This jet is $4m$-determined
since $x_3^2$ is equal to $t^{4m+2}$ for the singular curve and to zero for 
the regular ones. This jet is not enough, so we have to move to higher jets. 
A complete transversal in $J^{2n+1}$ is trivial for any $n$. Consider it 
in $J^{2n}\quad (n\le 2m)$. It is equivalent to
$$(G_2, (t^2, ct^{2n}, t^{2m+1}+dt^{2n}, et^{2n})).\eqno (12)$$
If $e\ne 0$ we obtain the normal form $4.2.4$. If $e=0$ and $c\ne 0\ne d$
we obtain the normal form $4.2.5$ if $n<2m$. If $n=2m$ we can eliminate
$t^{2n}$ in the third co-ordinate, since the tangent space to the
$A^{(2n)}$-orbit contains the vectors $(2t^{2m+1}, 0, (2m+1)t^{4n})$
and $(t^{2m+1}, 0, 0)$. Hence we obtain the normal form $4.2.6$ $(n=2m)$.
If $d=0$ and $c\ne 0$ we also obtain the normal form $4.2.6$. Now suppose
$c=0$ and $d\ne 0$. We obtain the $(2n)$-jet $(t^2, 0, t^{2m+1}+t^{2n})$, where
$n<2m$. We have to consider a complete transversal in $J^{2l}\quad (n<l\le 2m)$.
It is equivalent to $(t^2, pt^{2l}, t^{2m+1}+t^{2n}, qt^{2l})$ since
the tangent space to the $A^{(2l)}$-orbit contains 
$(2t^{2(1+l-n)}, (2m+1)t^{2(m+l-n)+1}+2nt^{2l})$ and 
$(0, t^{2(m+l-n)+1}+t^{2l})$. If $q\ne 0$ we obtain the normal
form $4.2.7$. If $q=0$ and $p\ne 0$ we obtain $4.2.8$. If $p=0$ and $q=0$, for 
every $l$ such that $n<l\le 2m$, we obtain $4.2.9$. If $c=d=e=0$ in (12) we obtain 
the normal form $4.2.10$.

Now suppose $b=0$ and $a\ne 0$ in (11). We have the $(2m+1)$-jet
$(t^2, t^{2m+1})$. It is $(4m+2)$-determined, since $x_1 x_2^2$ is equal to 
$t^{4m+4}$ for the singular curve and to zero for regular ones. A complete 
transversal in $J^{2n+1}$ is trivial for any $n$. In $J^{2n}$, where $m<n\le 2m+1$, 
it is equivalent to  $$G_2, (t^2, t^{2m+1}+ct^{2n}, dt^{2n}). \eqno (13)$$ 
If $d\ne 0$ we obtain the normal form $4.2.11$. If $d=0$ and $c\ne 0$ we move
to higher jets. We suppose $n<2m+1$, since if $n=2m+1$ the tangent space
to the $A^{(4m+2)}$-orbit contains the vectors: $(2t^{2m+3}, (2m+1)t^{4m+2})$ 
and $(t^{2m+3}, 0)$. A complete transversal in $J^{2l}$, where 
$n<l\le 2m+1$, is equivalent to $(t^2, t^{2m+1}+t^{2n}, et^{2l})$. If $e\ne 0$
we obtain the normal form $4.2.12$. If $e=0$ for every $l$ such that $n<l\le 2m+1$,
we obtain the normal form $4.2.13$. If in (13) $c=d=0$ for every $n$
such that $m<n\le 2m$, we obtain the normal form $4.2.14$. 

Now we shall consider a complete transversal in $J^{2m}$. It is equivalent to
$$(G_2, (t^2, at^{2m}, bt^{2m})).\eqno (14)$$ 
Suppose $a\ne 0$ or $b\ne 0$.
This jet is not finite determined, so we need to move to higher jets.
In $J^{2n}\quad (m<n)$ a complete transversal is trivial for any $n$. 

First, suppose $b\ne 0$. Then this $(2m)$-jet is equivalent to
$(G_2, (t^2, 0, t^{2m}))$. In $J^{2n+1}\quad (m\le n)$ a complete transversal
is equivalent to 
$$(t^2, ct^{2n+1}, t^{2m}+dt^{2n+1}, et^{2n+1}).$$
If $e\ne 0$ we obtain the normal form $4.2.15$, which is $(2n+1)$-determined.
Now suppose $e=0$. If $c\ne 0$ and $d\ne 0$ we obtain the normal form $4.2.18$,
where $l=n$. If $c\ne 0$ and $d=0$ we obtain the normal form $4.2.16$. These
both forms are $(2n+1)$-determined. Now suppose $c=e=0$ and $d\ne 0$. 
Then we have the $(2n+1)$-jet $(t^2, 0, t^{2m}+t^{2n+1})$. This jet
is $(2(m+n)-1)$-determined, since $x_3^2-x_1^m x_3$ is equal to $t^{2(m+n)+1}$ 
for the singular component and to zero for the regular ones. A complete transversal in
$J^{2l}$ is trivial for any $l$. In $J^{2l+1}$ $(n<l<m+n)$ it is equivalent to
$$G_2, (t^2, pt^{2l+1}, t^{2m}+t^{2n+1}, qt^{2l+1})$$
since the tangent space to the $A^{(2l+1)}$-orbit contains the vectors
$(2t^{2(1+l-n)}, 0, 2mt^{2(m+l-n)}+(2n+1)t^{2l+1})$ and
$(0, 0, t^{2(m+l-n)}+t^{2l+1})$.
If $q\ne 0$ we obtain the normal form $4.2.17$. If $q=0$ and $p\ne 0$ we
obtain the normal form $4.2.18$. If $p=q=0$ for every $l$ such that
$n<l<m+n$, we obtain the normal form $4.2.19$. 

Now, suppose $b=0$ and $a\ne 0$ in (14). We obtain the $2m$-jet 
$(G_2, (t^2, t^{2m}))$. A complete transversal in $J^{2n}$ is trivial for any $n$.
In $J^{2n+1}$ $m\le n$ it is equivalent to 
$(G_2, (t^2, t^{2m}+ct^{2n+1}, dt^{2n+1})).$
Suppose $d\ne 0$. We obtain the normal form $4.2.20$ which 
is $(2n+1)$-determined. If $d=0$ and $c\ne 0$
we have the $(2n+1)$-jet $(t^2, t^{2m}+t^{2n+1})$ which is
$(2(m+n)+1)$-determined since $x_1 x_2^2$ is equal to $t^{2(m+n)+3}$ for the singular
component and to zero for the regular ones. A complete transversal 
in $J^{2l+1}$ ($n<l\le m+n$) is equivalent to $(t^2, t^{2m}+t^{2n+1}, et^{2l+1}).$
If $e\ne 0$ we obtain the normal form $4.2.21$. If $e=0$ for each $l$ such that
$n<l\le m+n$, we obtain the normal form $4.2.22$.
\end{proof}

Now we shall consider multigerms with the regular part $G_n$ and the
multiplicity of the singular component is greater than or equals 3.
First we shall consider the 3-jet $(G_n, (\underbrace{0,\ldots, 0}_{n},t^3))$.

\begin{lemma}\label{Lem3_2.6}
Any multigerm with the 5-jet ($n>1$)
$$(G_n, (a_1t^3+b_1t^4+c_1t^5,\ldots, a_{n+1}t^3+b_{n+1}t^4+c_{n+1}t^5))$$
is not simple.
\end{lemma}

\begin{proof}
The dimension of the submanifold of such jets in $J^5$ is equal to $3n+3$. 
Let us estimate the dimension of the orbit at points
of the submanifold under the action of the stabilisator of the regular part. 
It is generated by the images of the following vectors from $L^{(5)}$:
$(\underbrace{0,\ldots, 0}_{k-1}, x_k,0,\ldots, 0), 1\le k\le n+1,$ and
$(\underbrace{0,\ldots, 0}_{k-1}, x_{n+1},0,\ldots, 0), 1\le k\le n$.
Hence the dimension of the orbit under the action of elements from 
$L^{(5)}$ is less than or equals $(n+1)+n=2n+1$. 
Only $t$, $t^2$ and $t^3$ can give rise to non-zero
vectors. Hence the dimension of the tangent space to the $R^{(5)}$-orbit 
is less than or equals 3. Now one can see that the dimension of the
tangent space to the $A^{(5)}$-orbit at points of the submanifold is
less than or equals $2n+4$. By similar resoning as in Lemma~\ref{Lem3_2.3}
we conclude that if the multigerm is simple then $3n+3\le 2n+4$, i.e.,
$n\le 1$. We have a contradiction with the condition of the Lemma.
\end{proof}

\begin{lemma}
Any simple multigerm with the 3-jet 
$(G_n, (\underbrace{0,\ldots, 0}_{n},t^3)),\quad n>1,$
is equivalent to one of $4.1.6$--$4.1.10$.
\end{lemma}

\begin{proof}
A complete transversal in $J^4$ is equivalent to 
$(a_1t^4,\ldots, a_nt^4, t^3, bt^4)$
If $b\ne 0$ we have the 4-jet $(0,\ldots, 0, t^3, t^4)$.
A complete transversal in $J^5$ is equivalent to
$(c_1t^5,\ldots, c_nt^5, t^3, t^4,dt^5).$
For any $c$ and $d$ this jet is 5-determined. If $d\ne 0$ we obtain the 
normal form $4.1.6$. If $d=0$ we obtain the normal form $4.1.7$. 
If $b=0$ the 4-jet is equivalent to
$(\underbrace{t^4,\ldots, t^4}_{k},\underbrace{0,\ldots, 0}_{n-k}, t^3)$.
Using Lemma~\ref{Lem3_2.6} we note that there is only one possibility
for the 5-jet:
$$(\underbrace{t^4,\ldots, t^4}_{k},\underbrace{0,\ldots, 0}_{n-k},t^3,
                                                   dt^5),\quad d\ne 0.$$
If $k\ne 0$ this jet is 5-determined and we obtain the normal form $4.1.8$.
If $k=0$ we move to higher jets. A complete transversal in $J^6$ is trivial.
A complete transversal in $J^7$ is equivalent to
$(p_1t^7,\ldots, p_nt^7, t^3, t^5, qt^7).$
If $q\ne 0$ this 7-jet is 7-determined and we obtain the normal form $4.1.9$.
If $q=0$ and there exists $i$ such that $p_i\ne 0$ we obtain the normal 
form $4.1.10$.
If $q=0$, and for each $i$ $p_i=0$ we obtain the normal form $4.1.8$ with $k=0$
which is 7-determined.
\end{proof}

Now we shall complete the classification of multigerms with 
non-zero 3-jet of the singular component and the regular part $G_n$. 
At first, note that $n=2$ since the 3-jet satisfies the conditions of 
Lemma~\ref{Lem3_2.3}. It is obvious, that we have the following
possibilities
for 3-jet of the singular component: $(0, 0, t^3)$, $(t^3, t^3, 0)$ and
$(t^3, 0, 0)$. $(0, t^3, 0)$ we don' t consider cause it can be obtained
from the last 3-jet by permutations of co-ordinates and components.

\begin{lemma}
Any simple multigerm with the 3-jet as above is equvalent to one of
$4.3.1$--$4.3.4$.
\end{lemma}

\begin{proof}
Using Lemma~\ref{Lem3_2.6} we conclude that 5-jet is equivalent to
$(a_1t^3, a_2t^3, t^4, t^5)$. If $a_1\ne 0$ and $a_2\ne 0$ then we obtain
the normal form $4.3.1$ which is 5-determined. Let $a_2=0$ and $a_1\ne 0$.
We have the 5-jet $(t^3, 0, t^4, t^5)$. This jet is 6-determined, so we
have to consider a complete transversal in $J^6$. It is equivalent to
$(t^3, bt^6, t^4, t^5, ct^6)$
since the tangent space to the $A^{(6)}$-orbit contains the vectors:
$(3t^4, 0, 4t^5, 5t^6, 0)$, $(t^4, 0, 0, 0, 0)$, $(0, 0, t^5, 0, 0)$ and
$(3t^5, 0, 4t^6, 0)$.
If $c\ne 0$ in (10) we obtain the normal form $4.3.2$. If $c=0$ and $b\ne 0$ we 
obtain the normal form $4.3.3$. If $b=c=0$ we obtain the normal form $4.3.4$. 
\end{proof}

Now we shall consider multigerms with zero 3-jet of the singular component.
Using Lemma~\ref{Lem3_2.3} we conclude that $n=2$, moreover by
Lemma~\ref{Lem3_2.6} there is only one
possibility for the 5-jet: $(G_2, (0, 0, t^4, t^5))$.

\begin{lemma}\label{Lem3_2.9}
Any multigerm with the 7-jet
$$(G_2, (a_1t^4+b_1t^5+c_1t^6+d_1t^7,\ldots, a_4t^4+b_4t^5+c_4t^6+d_4t^7))$$
is not simple.   
\end{lemma}

\begin{proof}
The dimension of the submanifold of such jets in $J^7$ is equal to $16$. 
Let us estimate the dimension of the orbit at points
of the submanifold under the action of the stabilisator of the regular part. 
It is generated by the images of the vectors:
$(\underbrace{0,\ldots, 0}_{i-1},x_i)\quad 1\le i\le 4$, $(0, 0, 0, x_3)$,
$(0, 0, x_4, 0)$; $(x_j, 0, 0, 0)$, $(0, x_j, 0, 0),\quad j=1,2$. 
Hence the dimension of the orbit under the action of elements from $L^{(7)}$ 
is less than or equals $4+3+3=10$. Only $t$, $t^2$, $t^3$ and $t^4$ can 
give rise to non-zero vectors. Hence the dimension of the tangent space to 
the orbit under the action of elements from $R^{(7)}$ is
less than or equals 4. Now one can see that the dimension of the
tangent space to the orbit at points of 
the submanifold is less than or equals $10+4=14$. Since $16\ge 14$, the multigerm
is not simple.
\end{proof}

\begin{lemma}
Any simple multigerm with the 5-jet $(G_2, (0, 0, t^4, t^5))$
is equivalent to one of $4.3.5$--$4.3.11$.
\end{lemma}

\begin{proof}
A complete transversal in $J^6$ is equivalent to 
$$(at^6, bt^6, t^4, t^5+dt^6, ct^6).\eqno (15)$$
Since the tangent space to the $A^{(6)}$-orbit contains the vectors 
$(0, 0, 4t^5, 5t^6)$, $(0, 0, t^5+dt^6, 0)$ and $(0, 0, 4t^6, 0)$, 
we can suppose $d=0$. If $c\ne 0$ we have the 6-jet 
$(0, 0, t^4, t^5, t^6)$.
This jet is 7-determined, so we move to $J^7$. A complete transversal
is equivalent to $(pt^7, qt^7, t^4, t^5, t^6, rt^7)$. If $r\ne 0$ we obtain
the normal form $4.3.5$. If $r=0$ but $p\ne 0$ and $q\ne 0$ we obtain 
the normal form $4.3.6$. If $p=r=0$ and $q\ne 0$ or $q=r=0$ and $p\ne 0$ 
we obtain the normal form $4.3.7$ (these both cases are equivalent since we 
can permutate co-ordinates and regular components). 
If $p=q=r=0$ we obtain the normal form $4.3.8$.

If $c=0$ in (15), we have (by similar reasoning) three possibilities 
for the 6-jet: $(t^6, t^6, t^4, t^5)$, $(t^6, 0, t^4, t^5)$ and $(0, 0, t^4, t^5)$.
Using Lemma~\ref{Lem3_2.9} we see that there is only one multigerm
for each of these three 6-jets~--- the normal forms $4.3.9$--$4.3.11$. 
All of them are 7-determined. 
\end{proof}

We have completely considered simple multigerms with one singular component
and the regular part $G_n$.
Now we have to consider multigerms with one singular component and the 
regular part $$((t_1, 0), (t_2, t_2^m)).\eqno (16)$$

Lemma~\ref{Lem3_0.2}
states that a multigerm with the 2-jet
$$((t_1, 0), (t_2, t_2^2), (t_3, \alpha t_3^2))\eqno (17)$$
is not simple. Since the third component of our multigerm is singular
we obtain $m=2$ in (16)  (if $m>2$ the 2-jet of our multigerm is 
adjacent to $((t_1, 0), (t_2, 0), (0, at^2))$ which is adjacent
to the family (17) ). Moreover there is only one possibility for the 2-jet:
$$((t_1, 0, 0), (t_2, t_2^2, 0), (0, 0, t^2)).\eqno (18)$$

\begin{lemma}
Any simple multigerm with the 2-jet $(18)$ is equivalent
to one of $4.4.1$--$4.4.3$.
\end{lemma}

\begin{proof}
The 2-jet (18) is not finite determined. 
A complete transversal in $J^{2m}$ is trivial for any $m$. In $J^{2m+1}$ it 
is equivalent to 
$((t_1, 0), (t_2, t_2^2), (at^{2m+1}, bt^{2m+1}, t^2, ct^{2m+1})).$
If $c\ne 0$ we obtain the normal form $4.4.1$, which is $(2m+1)$-determined.
If $c=0$ and $b\ne 0$ then the jet  is equivalent to
$$(t_1, 0), (t_2+at_2^2, t_2^2), (0, t^{2m+1},  t^2).$$
Since the tangent space to the $A^{(2m+1)}$-orbit contains the
vectors (we write only non-zero (the second) component):
$(t_2^2+2at_2^3, 2t_2^3)$, $(t_2^3+at_2^4, 0)$, $(0, t_2^3+at_2^4)$,
$(t_2^4, 0)$ and $(0, t_2^4)$
then the tangent space also contains $(t_2^2, 0)$. Hence by the Mather Lemma 
we conclude that the $(2m+1)$-jet is equivalent to the normal form $4.4.2$.
If $b=c=0$ and $a\ne 0$ we obtain the normal form $4.4.3$.  
\end{proof}

\section{Multigerms with two singular components}

Lemma~\ref{Lem3_0.2} states that the family of 2-jets
$$((t, 0), (t_1, t_1^2), (t_2, \alpha t_2^2)) \eqno (19)$$
is not simple.

Suppose that the multigerm contains $n$ regular components and
two singular ones. Then the 1-jet of the multigerm is adjacent to
$(n+2)$ lines in ${\mathbb C}^n$. In the proof of Lemma~\ref{Lem3_1.1}
we saw that if $n>1$ then
this 1-jet is not simple. Hence the multigerm contains three components:
one regular and two singular. 

The 2-jet of the multigerm is equivalent to
$$((t, 0, 0), (0, t_1^2, 0), (0, 0, t_2^2)) \eqno (20) $$ 
since it can not be adjacent to (19).

\begin{lemma}\label{Lem3_3.1}
The multigerm with the 3-jet
$$((t, 0, 0), (a_1 t_1^2+b_1 t_1^3, a_2 t_1^2+b_2 t_1^3, a_3 t_1^2+ b_3 t_1^3),
(a_4 t_2^2+b_4 t_2^3, a_5 t_2^2+b_5 t_2^3, a_6 t_2^2+ b_6 t_2^3))$$
is not simple.
\end{lemma}

\begin{proof}
The dimension of the submanifold of such jets in $J^3$ equals 12.
Let us estimate the dimension of the orbit at points
of the submanifold under the action of the stabilisator of the regular part. 
The dimension of the orbit under the action of elements from $L^{(3)}$ 
is less than or equals $3+2+2=7$. The dimension of the orbit under the action
of elements from $R^{(3)}$ is less than or equals $2+2=4$
Now one can see that the dimension of the tangent space to the orbit 
at points of the submanifold is less than or equals $7+4 = 11 < 12$.
\end{proof}

We shall not write the regular component (we suppose it equals $(t, 0, 0)$).
Using Lemma~\ref{Lem3_3.1} we can suppose that the second singular
component equals $(0, 0, t_2^2, t_2^3)$.

\begin{lemma}
Any simple multigerm with the 2-jet $(20)$ is equvalent to
one of $5.1$--$5.6$. 
\end{lemma}

\begin{proof}
A complete transversal in $J^{2m}$ is trivial for any $m$. 
In $J^{2m+1}$ a complete transversal is equivalent to 
$$(at_1^{2m+1}, t_1^2, bt_1^{2m+1}, ct_1^{2m+1}, dt_1^{2m+1})$$
(we write only the second component, the first and the third components 
are not changed).
If $d\ne 0$ we obtain the normal form $5.1$. If  $d=0$ and $ c\ne 0$
we can suppose $c=1$.
The tangent space to the $A^{(2m+1)}$-orbit contains the vectors 
(we write only the singular components): 
$((0, 0, t_1^{2m+1}, 0), (0, 0, t_2^3, 0))$,
$((0,0,0,0), (0,0, 2t_2^3, 3t_2^4))$ and $((0,0,0,0), (0,0,0, t_2^4))$. 
So, using the Mather Lemma we can suppose $b=0$. If $a\ne 0$ we obtain 
the normal form $5.2$. If $a=0$ we obtain the normal form $5.3$.

Now suppose $c=d=0$. If $a\ne 0$ and $b\ne 0$ we obtain the normal
form $5.4$. If $b\ne 0$ and $a=0$ we obtain the normal form $5.5$. If $b=0$
and $a\ne 0$ we obtain the normal form $5.6$. All these normal
forms are $(2m+1)$-determined.
\end{proof}

Moscow State University, \\
Faculty of Mechanics and Mathematics.\\
Moscow, 119899 \\
Russia.\\
E-mail: kolgush@mccme.ru, sadykov@mccme.ru  
\end{document}